\documentstyle[12pt]{article}
\newtheorem{theorem}{Theorem}[section]
\newtheorem{prop}{Proposition}[section]

\newtheorem{lemma}{Lemma}[section]
\newtheorem{definition}{Definition}[section]

\newcommand{\U}{U_q(C_n^{(1)})}
\renewcommand{\a}{\alpha}

\newcommand{\la}{\lambda}
\newcommand{\be}{\beta}

\newcommand{\ep}{\epsilon}
\newcommand{\vep}{\varepsilon}

\title{Vertex Operators of Admissible Modules of $U_q(C^{(1)}_n)$}

\author{Naihuan Jing and Yoshitaka Koyama}
\begin{document}           
\date{October 10, 1996}
\maketitle                 
\begin{center}
Department of Mathematics\\
North Carolina State University\\
Raleigh, NC 27695-8205, USA
\end{center}

\begin{abstract}
Using our recent bosonic realization of $U_q(\widehat{sp}_{2n})$,
we construct explicitly the vertex operators for the level $-1/2$ modules of
$U_q(\widehat{sp}_{2n})$ using bosonic fields. Our method contains
a detailed analysis of all the $q$-intertwining relations.
\end{abstract}

\section{Introduction}
The $q$-vertex operators play a crucial role in the mathematical
formulation of the solvable lattice models (cf. \cite{kn:JbM})
and in the theory of $q$-KZ equations \cite{kn:FR}. In particular the
method works well when the explicit bosonization is available (cf.
\cite{kn:JMMN}).

The program of realization of vertex operators was carried out for many
quantum affine algebras of classical types for both untwisted and twisted
types at level one irreducible modules
in \cite{kn:K1} \cite{kn:JM1} \cite{kn:JKK} \cite{kn:JM2},
and level two irreducible modules of $U_q(\widehat{sl}_2)$ in \cite{kn:I}
\cite{kn:JM1}.

Recently we have given an explicit bosonic realization of the
quantum affine algebra $U_q(\widehat{sp}_{2n})$ by admissible
representations of level $-1/2$ in \cite{kn:JKM}. This realization is
a $q$-analogue of the Feingold-Frenkel construction \cite{kn:FF} of
affine symplectic algebras and also
generalizes the result of $U_q(\widehat{sl}_2)$ in \cite{kn:K2}.

In this paper we use this explicit realization to construct $q$-vertex
operators associated with the admissible representations. Admissible
representations are generalizations of the integrable (integral levels)
modules to those with levels of rational numbers \cite{kn:KW}.
Their characters are also modular functions.

Our construction has a new feature different from previous cases of
positive levels. The intertwining
relations are much harder to prove in our case due to appearance of poles
in the contraction functions, while in the level one case
they are consequences of some coproduct relations of Drinfeld generators.
We give a detail analysis of the dependency among the intertwining relations
and single out
the key relations between the $0$-vertex and the highest component of
the vertex operators. These key relations are then proved by a long chain of
$q$-commutation relations of the algebra and the operators. The method
clearly works for all other cases we mentioned above, and thus is completely
general by its own means.

It is still unclear what physical models our vertex operators will
describe. We only have an imbedding of
$V(\mu_i)$ into $V(\mu_j)\otimes V_z$ instead of isomorphism in level one
cases. It will be interesting to find out exactly
how the level $-1/2$-modules are imbedded in $\otimes^{\infty/2}V$.
It is clear that our $q$-vertex operator components generate a $q$-Weyl
algebra. It would also be interesting to construct $U_q(C_n^{(1)})$ from this
$q$-Weyl algebra, which is regarded as one direction of the so-called
boson-fermion correspondence (cf. \cite{kn:F}).
In this aspect the work of \cite{kn:D1}\cite{kn:D2}
may be useful, where the similar problem for $q$-Clifford algebras was
studied.

The paper is organized as follows. In section 2 we review the
Feingold-Frenkel construction and discuss the equivalent form in terms of
scaler fields. It is the latter form we quantized in \cite{kn:JKM} and recalled
in section 3. The main results are also stated in section 3.
In section 4 we start our analysis of the intertwining relations and reduce
them to a few key relations, and point out that we only need to
prove one of them. In the final section we prove the key
relation.

\section{Feingold-Frenkel construction of
	 $\widehat{sp}_{2n}$ and vetex operators.}

   In this section,
   we review the Feingold-Frenkel construction of $\widehat{sp}_{2n}$
   in terms of the bosonic $\beta$-$\gamma$ system
   and describe the relations between the bosonic $\beta$-$\gamma$ system
   and vertex operators associated with the vector representation.

 \subsection{Affine Lie algebras $\widehat{sp}_{2n}$.}

   We recall the affine Lie algebra $\widehat{sp}_{2n}$
   as follows.
   Most notations concerning affine Lie algebras follow \cite{kn:K}.
   Let $V$ be a $2n$-dimensional {\bf C}-vector space
   with basis $v_1, \cdots, v_n$,
   $v_{\overline 1}, \cdots, v_{\overline n}$.
   The Lie algebra $sp_{2n}$ is defined as
   the Lie subalgebra of $gl(V)$ generated by the elements:
   $$ E_i=E_{i i+1}-E_{\overline{i+1} \, \overline{i}}, \quad
      F_i=E_{i+1 i}-E_{\overline{i} \, \overline{i+1}}, \quad
      E_n=E_{n \overline{n}}, \quad
      F_n=E_{\overline{n} n}, $$
   $$ H_i=E_{ii}-E_{i+1i+1}
	  +E_{\overline{i+1} \, \overline{i+1}}
	  -E_{\overline{i} \, \overline{i}} , \quad
      H_n=E_{nn}-E_{\overline{n} \, \overline{n}} . $$
   where $E_{ij}$, $E_{\overline{i} \overline j}, \cdots$
   is the matrix unit of End$(V)$
   such that $E_{ij} v_k = \delta_{jk} v_i$, etc.
   We will refer to $V$ as the vector representation of $sp_{2n}$.
   As a vector space $\widehat{sp}_{2n}$ is given by
   $$ \widehat{sp}_{2n}=sp_{2n} \otimes {\bf C}[t,t^{-1}]
			\oplus {\bf C}K \oplus {\bf C}d. $$
   The Lie algebra structure is defined by
   $$ [X \otimes t^l, Y \otimes t^m]
     =[X,Y] \otimes t^{l+m} +l \delta_{l+m,0} K tr(XY) , \quad
      [X \otimes t^m, K]=0, $$
   $$ [d, K]=0, \quad
      [d, X \otimes t^m]=m X \otimes t^m \quad
      \mbox{for} \; X,Y \in sp_{2n} . $$
   We set the following generating function for $X \in sp_{2n}$.
   $$ X(z)=\sum_{m \in {\bf Z}} (X \otimes t^m) z^{-m} . $$

 \subsection{Feingold-Frenkel construction.}

   As an index set,
   we set $I={\bf Z}$ or ${\bf Z}+\frac{1}{2}$.
   We introduce $n$ copies of the bosonic $\beta$-$\gamma$ systems,
   which are the set of operators
   $\beta_i(k)$, $\gamma_i(k)$ $(i=1, \cdots, n, \; k \in I)$
   satisfying the defining relations:
   $$ [\beta_i(l) ,\gamma_j(m)]=\delta_{ij}\delta_{l+m,0}, \quad
      [\beta_i(l) ,\beta_j(m) ]=0 , \quad
      [\gamma_i(l),\gamma_j(m)]=0 , $$
   $$ \mbox{for} \quad 1 \leq i,j \leq n, \; l,m \in I . $$
   The Fock representation ${\cal F}_{\beta\mbox{-}\gamma}$
   is generated by the vacuum vector $|vac \rangle$ with
   the following defining relations.
   $$ \beta_i(l) |vac \rangle =0, \quad
      \gamma_i(m)|vac \rangle =0  \quad
      \mbox{for} \quad i=1, \cdots, n, \; l \geq 0 , \; m > 0 . $$
   Let $\beta_i(z)$, $\gamma_i(z)$ be the generating functions:
   $$ \beta_i(z) =\sum_{k \in I} \beta_i(k) z^{-k}, \;
      \gamma_i(z)=\sum_{k \in I} \gamma_i(k) z^{-k}. $$
   Following \cite{kn:FF},
   we define the action of
   $\widehat{sp}_{2n}$ on ${\cal F}_{\beta\mbox{-}\gamma}$ by
   $$ H_i(z) =:\beta_i(z)\gamma_i(z):-:\beta_{i+1}(z)\gamma_{i+1}(z): , $$
   $$ E_i(z) = \beta_i(z)     \gamma_{i+1}(z) , \quad
      F_i(z) = \beta_{i+1}(z) \gamma_i(z)     , \quad
      \mbox{for} \; i=1, \cdots, n-1, $$
   $$ H_n(z) =-:\beta_n(z)\gamma_n(z): , \quad
      E_n(z) =-\frac{1}{2}\beta_n(z) \beta_n(z) , \quad
      F_n(z) =\frac{1}{2} \gamma_n(z)\gamma_n(z) , $$
   $$ K=-\frac{1}{2}, \quad
      d=-\sum^n_{i=1} \sum_{k \in I} k :\beta_i(k) \gamma_i(-k): . $$
   where $: \quad :$ is the normal ordering defined by
   $$ :\beta_i(l) \gamma_i(m):
	  = \left \{ \begin{array}{ll}
		     \beta_i(l) \gamma_i(m)  & (l<m) \\
		     \frac{1}{2}( \beta_i(l) \gamma_i(m)
				 +\gamma_i(m) \beta_i(l))
		     & (l=m) \\
		     \gamma_i(l) \beta_i(m)  & (l>m)
		     \end{array} \right. $$
   Then we have the following theorem.

   \begin{theorem}(\cite{kn:FF})
   ${\cal F}_{\beta \mbox{-} \gamma}$ is a $\widehat{sp}_{2n}$-module
   and decomposed into
   $$ L(-\frac{1}{2}\Lambda_0)
      \oplus
      L(-\frac{3}{2}\Lambda_0+\Lambda_1-\frac{1}{2}\delta)
      \quad \mbox{for} \quad I={\bf Z}+\frac{1}{2}, $$
   $$ L(-\frac{1}{2}\Lambda_n)
      \oplus
      L(-\frac{3}{2}\Lambda_n+\Lambda_{n-1})
      \quad \mbox{for} \quad I={\bf Z}. $$
   Here $L(\lambda)$ is
   the irreducible highest weight $\widehat{sp}_{2n}$-module
   with the highest weight $\lambda$,
   $\Lambda_i$ is the $i$-th fundamental weight
   and $\delta$ is the canonical imaginary root.
   Highest weight vectors are given by
   $|vac\rangle$, $\be_1(-\frac{1}{2})|vac\rangle$,
   $|vac\rangle$, $\gamma_n(0)|vac\rangle$
   for
   $L(-\frac{1}{2}\Lambda_0)$,
   $L(-\frac{3}{2}\Lambda_0+\Lambda_1-\frac{1}{2}\delta)$,
   $L(-\frac{1}{2}\Lambda_n)$,
   $L(-\frac{3}{2}\Lambda_n+\Lambda_{n-1})$ respectively.
   \end{theorem}

   From now on we set $\mu_1=-\frac{1}{2}\Lambda_0$,
   $\mu_2=-\frac{3}{2}\Lambda_0+\Lambda_1-\frac{1}{2}\delta$,
   $\mu_3=-\frac{1}{2}\Lambda_n$, and
   $\mu_4=-\frac{3}{2}\Lambda_n+\Lambda_{n-1}$.

 \subsection{Vertex operators.}

   We first define $V_z$, the affinization of the vector representation V,
   which is a level $0$ representation of $\widehat{sp}_{2n}$.
   As a vector space
   $$ V_z=V \otimes {\bf C}[z,z^{-1}] . $$
   The action of $\widehat{sp}_{2n}$ is given by:
   $$ X(m).v \otimes z^n = X.v \otimes z^{m+n} , \quad
      d.v \otimes z^n = n v \otimes z^n        , \quad
      K.v \otimes z^n = 0 . $$
   We define the operators $\phi^V(z)$ as follows.
   $$ \phi^V(z)= \sum^n_{i=1} \gamma_i(z) \otimes v_i
		+\sum^n_{i=1} \beta_i(z) \otimes v_{\overline i}, $$
   which can be regarded as an operator
   from $L(\mu)$ to a suitable completion of $L(\mu') \otimes V_z$,
   where $(\mu,\mu')=(-\frac{1}{2}\Lambda_0,
   -\frac{3}{2}\Lambda_0+\Lambda_1-\frac{1}{2}\delta)$,
   $(-\frac{3}{2}\Lambda_0+\Lambda_1-\frac{1}{2}\delta,
     -\frac{1}{2}\Lambda_0)$,
   $(-\frac{1}{2}\Lambda_n, -\frac{3}{2}\Lambda_n+\Lambda_{n-1})$
   or
   $(-\frac{3}{2}\Lambda_n+\Lambda_{n-1}, -\frac{1}{2}\Lambda_n)$.
   Then we have the following theorem.

   \begin{theorem}
   $\phi^V(z)$ is intertwining operator,
   or the vertex operators.(cf. \cite{kn:TK}, \cite{kn:FR})
   \end{theorem}
   {\it Proof.}
   We can check the intertwining relations
   between $\phi^V(z)$ and the $\widehat{sp}_{2n}$ currents immediately.
   \hspace{\fill} $\Box$

 \subsection{Representations in terms of scaler fields.}

  In this subsection, we explain the $\beta$-$\gamma$ system
  in terms of bosonic scaler fields (\cite{kn:FMS}).
  Using this picture, we can rewrite the Feingold-Frenkel construction
  in terms of the bosonic scaler fields.
  The rewritten representation coincides precisely
  with the classical limit ($q=1$) of
  the representation to be constructed in section 3.
  Let $\phi_{i1}(z)$, $\phi_{i2}(z)$ $(i=1, \cdots, n)$ be
  independent bosonic scaler fields normalized by
  $$ \phi_{i1}(z) \phi_{i1}(w) \sim -\log(z-w) , $$
  $$ \phi_{i2}(z) \phi_{i2}(w) \sim  \log(z-w) . $$
  It is known that
  the bosonic $\beta$-$\gamma$ system can be
  expressed by $\phi_{i1}(z)$ and $\phi_{i2}(z)$ as follows (\cite{kn:FMS}):
  $$ \beta_i(z)
     \longrightarrow
     :\partial_z \phi_{i2}(z) e^{\phi_{i1}(z)+\phi_{i2}(z)}: \; , $$
  $$ \gamma_i(z)
     \longrightarrow
     :e^{-\phi_{i1}(z)-\phi_{i2}(z)}: \; . $$
  We then obtain the Fock space of the bosonic $\beta$-$\gamma$ system
  by taking Ker$\,Q_i^-$ for all $i=1,\cdots,n$ with some charge constraints.
  Here $Q_i^-$ is the operator $\oint :e^{-\phi_{i2}(z)}: dz$.
  The correspondence with $Y_i(z)$ and $Z_i(z)$ in section 3
  is as follows:
  $$ Y^{\pm}_a(z) \stackrel{q \rightarrow 1}{\longrightarrow}
     :e^{\pm 2 (\phi_{i1}(z)-\phi_{i2}(z))}: \; , \quad
     Z^{\pm}_i(z) \longrightarrow
     :e^{\pm   \phi_{i2}(z)}: \; . $$

\section{Bosonic construction of $U_q(C^{(1)}_n)$
	 and $q$-vertex operators.}

 In this section,
 we review the bosonic construction of $U_q(C^{(1)}_n)$ in \cite{kn:JKM}
 and give a bosonization formula of the corresponding $q$-vertex operators.

 \subsection{Quantum affine algebras $U_q(C^{(1)}_n)$}

  Let $\widehat{P}$ be a free ${\bf Z}$-lattice of rank $n+2$
  and we denote the basis by $\vep_1, \cdots, \vep_n, d, \delta$.
  The nondegenerate inner product on $\widehat{P}$ is given  by
   $$ (\vep_i|\vep_j)=\frac{1}{2}\delta_{ij}, \;
      (d|\delta)=1, \;
      (d|d)=(\delta|\delta)=(\vep_i|d)=(\vep_i|\delta)=0 . $$
  We set the simple roots $\a_i$ and coroots $h_i$ by
   $$ \a_0=\delta-2\vep_1, \;
      \a_i=\vep_i-\vep_{i+1} \; (i=1, \cdots, n-1) , \;
      \a_n=2\vep_n, $$
   $$  h_0=\delta-2\vep_1, \;
       h_i=2\vep_i-2\vep_{i+1} \; (i=1, \cdots, n-1) , \;
       h_n=2\vep_n. $$
  Then the matrix $( (h_i|\a_j) )_{0 \leq i,j \leq n}$
  is the generalized Cartan matrix of type $C^{(1)}_n$.
  The fundamental weights $\Lambda_0$, $\cdots$, $\Lambda_n$
  and their classical parts
  $\lambda_1$, $\cdots$, $\lambda_n$ are given by
  $\Lambda_0=d$, $\Lambda_1=d+\vep_1$, $\Lambda_2=d+\vep_1+\vep_2$,
  $\cdots$, $\Lambda_n=d+\vep_1+\vep_2+\cdots+\vep_n$
  and $\lambda_i=\Lambda_i-d$ $(i=1, \cdots, n)$.
  Let $\widehat{P}^*$, $P$ and $Q$ be
  the sublattices of $\widehat{P}$ defined as follows.
  \begin{equation}
  \begin{array}{rcl}
  \ \widehat{P}^* &=& {\bf Z} h_0 \oplus \cdots \oplus {\bf Z} h_n
		      \oplus {\bf Z} d , \\
  \ P &=& {\bf Z} \vep_1 \oplus \cdots \oplus {\bf Z} \vep_n
	 ={\bf Z} \lambda_1 \oplus \cdots \oplus {\bf Z} \lambda_n , \\
  \ Q &=& {\bf Z} \a_1   \oplus \cdots \oplus {\bf Z} \a_n   .
  \end{array}
  \end{equation}

   The quantum affine algebra $U_q(C_n^{(1)})$ is
   the associative algebra with 1 over ${\bf C}(q^{1/2})$
   generated by the elements
   $q^h$ $(h \in \widehat{P}^*)$, $e_i$, $f_i$, $(i=0, 1, \cdots, n)$
   satisfying the following defining relations.
   $$ q^h=1 \quad \hbox{for $h=0$,} $$
   $$ q^{h+h'}=q^h q^{h'} \quad \hbox{for $h,h'\in \widehat{P}^*$,} $$
   $$ q^h e_i q^{-h} = q^{(h | \alpha_i)} e_i
      \quad \hbox{and} \quad
      q^h f_i q^{-h} = q^{-(h | \alpha_i )}f_i, $$
   $$ [e_i,f_j]=\delta_{ij}{{t_i-t_i^{-1}}\over{q_i-q_i^{-1}}}, $$
   $$ \sum_{k=0}^{b}
      \frac{(-1)^k}{[k]_i![b-k]_i!} e_i^{k} e_j e_i^{b-k}=0,
      \quad
      \sum_{k=0}^{b}
      \frac{(-1)^k}{[k]_i![b-k]_i!} f_i^{k} f_j f_i^{b-k}=0
      \quad \mbox{for} \; i \neq j. $$
   Here
   $$ b=1-(h_i|\a_j), \quad
      [k]_i=\frac{q_i^k-q_i^{-k}}{q_i-q_i^{-1}}, \quad
      [k]_i!=[1]_i[2]_i \cdots [k]_i, \quad
      q_i=q^{{(\alpha_i,\alpha_i)\over2}} \; \hbox{ and } \;
      t_i=q_i^{h_i}=q^{\a_i}. $$
   In this paper, we use the following comultiplication.
   $$ \Delta(q^h)=q^h \otimes q^h, \;
      \Delta(e_i)=e_i \otimes 1 + t_i \otimes e_i, \;
      \Delta(e_i)=f_i \otimes t_i^{-1} + 1 \otimes f_i. $$
   Throughout this paper, we denote
   $U_q(C_n^{(1)})$ by $U_q$
   and denote by $V(\lambda)$ the irreducible highest weight
   $U_q$-module with highest weight $\lambda$.
   We fix a highest weight vector of $V(\lambda)$
   and denote it by $|\lambda \rangle$.
   If $W_i$ $(i=1,2)$ has a weight decomposition
   $W_i=\oplus_{\nu} W_{i,\nu}$,
   their completed tensor product is then defined by
   $$ W_1 \widehat{\otimes} W_2
      =\bigoplus_{\nu}
      (\prod_{\nu=\nu_1+\nu_2} W_{1,\nu_1} \otimes W_{2,\nu_2}) . $$

  \subsection{Bosonic construction of $U_q(C^{(1)}_n)$.}

   Let $a_i(m)$ and $b_i(m)$ be the operators
   satisfying the following defining relations.

   \begin{equation}
   \begin{array}{rcl}
   \ [a_i(m), a_j(l)]&=&\delta_{m+l, 0}
			\frac{[m(\a_i|\a_j)][-\frac{m}{2}]}{m}, \\
   \ [b_i(m), b_j(l)]&=&m\delta_{ij}\delta_{m+l,0},\\
   \ [a_i(m), b_j(l)]&=&0 .
   \end{array}
   \end{equation}
   We define the Fock space
   ${\cal F}^a_\a$ and ${\cal F}^{(b,i)}_\be$
   for $\a \in P+\frac{\bf Z}{2}\la_n$, $\be \in P$
   by the defining relations
   $$ a_i(m) | \a, \be \rangle = 0 \quad (m>0)
      \; , \quad
      b_i(m) | \a, \be \rangle = 0 \quad (m>0) \;, $$
   $$ a_i(0) | \a, \be \rangle = (\a_i|\a) | \a, \be \rangle
      \; , \quad
      b_i(0) | \a, \be \rangle = (2\vep_i|\be) | \a, \be \rangle \; , $$
   where $| \a, \be \rangle$ is the vacuum vector.
   The grading operator ${d}$ is defined by
   $$ d. | \alpha,\beta \rangle
      =( (\alpha|\alpha)-(\beta|\beta-\lambda_n))
	|\alpha,\beta \rangle . $$
   We set
   $$ \widetilde{\cal F}=
	  \bigoplus_{\a\in P+\frac12{\bf Z}\la_n,\be \in P}
	  {\cal F}_{\a, \be} . $$
   Let $e^{a_i}=e^{\a_i}$ and $e^{b_i}$ ($1\leq i\leq n$) be operators
   on ${\widetilde{\cal F}}$ given by:
   $$ e^{\vep_i}|\a,\be \rangle= |\a+\vep_i,\be \rangle
      \quad , \quad
      e^{b_i}|\a, \be \rangle= |\a,\be+\vep_i \rangle . $$
   Let $: \quad :$ be the usual bosonic normal ordering defined by
   $$ :a_i(m) a_j(l): = a_i(m) a_j(l) \, (m \leq l) ,
   \; a_j(l) a_i(m) \, (m>l) , $$
   $$ :e^{a} a_i(0):=:a_i(0) e^{a}:=e^{a} a_i(0) \, . $$
   and similar normal products for the $b_i(m)$'s.
   Let $\partial=\partial_{q^{1/2}}$ be the $q$-difference operator:
   $$ \partial_{q^{1/2}}(f(z))
   =\frac{f(q^{1/2}z)-f(q^{-1/2}z)}{(q^{1/2}-q^{-1/2})z} $$
   We introduce the following operators.
\begin{eqnarray*}
&Y_i^{\pm}(z)&=
	 \exp ( \pm \sum^{\infty}_{k=1}
		 \frac{a_i(k)}{[-\frac{1}{2}k]} q^{\pm \frac{k}{4}} z^k)
	 \exp ( \mp \sum^{\infty}_{k=1}
	   \frac{a_i(k)   }{[-\frac{1}{2}k]} q^{\pm \frac{k}{4}} z^{-k})
	 e^{\pm a_i} z^{\mp 2a_i(0)} ,\\
&Z^{\pm}_i(z)&=
	     \exp ( \pm \sum^{\infty}_{k=1} \frac{b_i(-k)}{k} z^k)
	     \exp ( \mp \sum^{\infty}_{k=1} \frac{b_i(k)}{k} z^{-k})
	     e^{\pm b_i} z^{\pm b_i(0)} .
\end{eqnarray*}
   We define the operator $x_i^{\pm}(m)$ $(i=1, \cdots, n, m \in {\bf Z})$
   by the following generating function
   $X_i^{\pm}(z)=\sum_{m \in {\bf Z}} x_i^{\pm}(m) z^{-m-1}$.
\begin{eqnarray*}
X_i^+(z) &=&
	     \partial Z^+_i(z)Z^-_{i+1}(z)Y_i^+(z), \qquad i=1, \cdots n-1\\
X_i^-(z) &=&
	     Z^-_i(z)\partial Z^+_{i+1}(z)Y^-_i(z), \qquad i=1, \cdots n-1\\
X_n^+(z) &=&
	     \left(
	     \frac{1}{q^{\frac{1}{2}}+ q^{-\frac{1}{2}}}
	     :Z^+_n(z) \partial^2_z Z^+_n(z):
     -:\partial_z Y^+_b(q^{\frac{1}{2}}z)
	      \partial_z Y^+_b(q^{-\frac{1}{2}}z):
	     \right)
	     Y^+_a(z) \\
X_n^-(z) &=&
	     \frac{1}
		  {q^{\frac{1}{2}}+ q^{-\frac{1}{2}}}
	     :Z^-_n(q^{\frac{1}{2}}z)Z^-_n(q^{-\frac{1}{2}}z): Y^-_n(z).
\end{eqnarray*}

\noindent{\it Remark}. Our $a_i(k)$ differs from that in \cite{kn:JKM}, where
we took $a_i(k)/[d_i]$ for $a_i(k)$.

  \begin{theorem} \label{T:1} (\cite{kn:JKM})
  $\widetilde{\cal F}$ is a $U_q$-module of level $-\frac{1}{2}$
  under the action defined by
   $$ t_i \longmapsto q^{\a_i} , \quad
      e_i \longmapsto x^+_i(0), \quad
      f_i \longmapsto x^-_i(0)  \quad \mbox{for} \; i=1,\cdots,n \; , $$
   $$ t_0 \longmapsto q^{-\frac{1}{2}} K_{\theta}^{-1}
      \quad (K_{\theta}=q^{2\a_1+\cdots+2\a_{n-1}+\a_n}), $$
   $$ e_0 \longmapsto
	  \frac{(-1)^n}{[2]_1}
	  [x_1^-(0), \cdots, x_n^-(0),
	   x_{n-1}^-(0), \cdots, x_2^-(0), x_1^-(1)
	  ]_{q^{-1/2}\cdots q^{-1/2}\, q^{-1}\, q^{-1/2}\cdots q^{-1/2}\, 1}
	  K_{\theta}^{-1}, $$
   $$ f_0 \longmapsto
	  \frac{(-q)^n}{[2]_1}
	  [x_1^+(0), \cdots, x_n^+(0),
	   x_{n-1}^+(0), \cdots, x_2^+(0), x_1^+(-1)
	  ]_{q^{-1/2}\cdots q^{-1/2}\, q^{-1}\, q^{-1/2}\cdots q^{-1/2}\, 1}
	  K_{\theta}. $$
  \end{theorem}
  Furthermore, we know that $\widetilde{\cal F}$ contains
  the four irreducible highest weight modules (\cite{kn:JKM}).
  Let ${\cal F}_{\a,\be}^1$ be the subspace of ${\cal F}_{\a,\be}$
  generated by $a_i(m)$ $(i=1,\cdots,n \hbox{ and } m \in {\bf Z})$.
  Simularly, let ${\cal F}_{\a,\be}^{2,j}$ $(j=1, \cdots, n)$ be
  the subspace of ${\cal F}_{\a,\be}$ generated by
  $b_j(m)$ $(m \in {\bf Z})$.
  We can define the following isomorphism by
  $|\a,\be\rangle \otimes |\a',\be' \rangle
   \rightarrow |\a+\a',\be+\be' \rangle$.
  $$ {\cal F}^1_{\a,0} \otimes
     {\cal F}^{2,1}_{0,\be_1} \otimes \cdots \otimes {\cal F}^{2,n}_{0,\be_n}
     \longrightarrow
     {\cal F}_{\a, \be_1+\cdots+\be_n} . $$
  Let $Q_j^-$ be the operator from ${\cal F}^{2,j}_{\a, \be}$
  to ${\cal F}^{2,j}_{\a, \be-\vep_j}$ defined by
  \[ Q_i^-=\frac{1}{2 \pi \sqrt{-1}} \oint Z^-_i(z) dz .\]
  We set subspaces ${\cal F}_i$ $(i=1,2,3,4)$
  of $\widetilde{\cal F}$ as follows.

\[\begin{array}{llll}
  \displaystyle{
 {\cal F}_1 = \bigoplus_{\alpha \in Q}
		   {\cal F}'_{\alpha,\alpha}},
    &&\displaystyle{{\cal F}_2 = \bigoplus_{\alpha \in Q}
		   {\cal F}'_{\alpha+\vep_1,\alpha+\vep_1}}\\
    \displaystyle{{\cal F}_3 = \bigoplus_{\alpha \in Q}
		   {\cal F}'_{\alpha-\frac{1}{2}\lambda_n, \alpha}},
   &&\displaystyle{{\cal F}_4 = \bigoplus_{\alpha \in Q +\vep_n}
		   {\cal F}'_{\alpha-\frac{1}{2}\lambda_n, \alpha}},
\end{array} \]
   where
   $$ {\cal F}'_{\alpha,\beta}
     ={\cal F}^{1}_{\alpha,0}
      \otimes
      \bigotimes^n_{j=1}
      {\mbox Ker}_{{\cal F}^{2,j}_{0, l_j \vep_j}} Q^-_j , $$
   for $ \beta=l_1 \vep_1 + \cdots + l_n \vep_n $.
   Then we have the following theorem.

   \begin{theorem}(\cite{kn:JKM})
   Each ${\cal F}_i$ $(i=1,2,3,4)$ is
   an irreducible highest weight $U_q$-module
   isomorphic to $V(\mu_i)$,
   The highest weight vectors are given by
   $ | \mu_1 \rangle
    =|0, 0 \rangle $,
   $ | \mu_2 \rangle
    = b_1(-1) | \lambda_1, \lambda_1 \rangle $,
   $ | \mu_3 \rangle
    =| -\frac{1}{2}\lambda_n, 0 \rangle $,
   $ | \mu_4 \rangle
    =|-\frac{1}{2}\lambda_n-\vep_n, -\vep_n \rangle $.
   \end{theorem}

 \subsection{$q$-vertex operators}

   We recall the evaluation $\U$-module $V_z$.
   Let $V$ be a $2n$ dimensional vector space
   and $v_1, \cdots, v_n$, $v_{\overline 1}, \cdots, v_{\overline n}$
   be basis of $V$.
   Set $V_z=V \otimes {\bf C}[z,z^{-1}]$.
   The action of $\U$ on $V_z$ is defined as follows.
   $$ e_0=E_{\overline{1} 1} \otimes z, \quad
      e_i=E_{i i+1}-E_{\overline{i+1} \, \overline{i}}
      \quad (\mbox{for} \; i=1, \cdots, n-1), \quad
      e_n=E_{n \overline{n}}, $$
   $$ f_0=E_{1 \overline{1}} \otimes z^{-1}, \quad
      f_i=E_{i+1 i}-E_{\overline{i} \, \overline{i+1}},
      \quad (\mbox{for} \; i=1, \cdots, n-1), \quad
      f_n=E_{\overline{n} n}, $$
   $$ q^h v \otimes z^n = q^{(h|wt(v)+n\delta)} v \otimes z^n \quad
      \mbox{for} \; v=v_1, \cdots, v_n,
		    v_{\overline 1}, \cdots, v_{\overline n},$$
   where $wt(v)$ is given by
   $ wt(v_i)=\vep_i , \quad wt(v_{\overline i})=-\vep_i $.

   \begin{definition} (\cite{kn:FR})
   The $q$-vertex operator is
   a $U_q$-homomorphism of one of the following types.

    Type  I :
    $$ {\Phi}^{\mu V}_{\lambda}(z):
	  V(\lambda) \longrightarrow V(\mu) \widehat{\otimes} V_z     $$

    Type II :
    $$ ^V{\Phi}^{\mu}_{\lambda}(z):
	  V(\lambda) \longrightarrow V_z \widehat{\otimes} V(\mu)     $$
   \end{definition}
   We define
   the components ${\Phi}^{\mu V}_{\lambda}{}_{i} (z)$,
   ${\Phi}^{\mu V}_{\lambda}{}_{\overline{i}} (z)$
   of the $q$-vertex operators as follows.
   $$ {\Phi}^{\mu V}_{\lambda}(z)
     = \sum_{i=1}^n
	{\Phi}^{\mu V}_{\lambda}{}_i (z) \otimes v_i
      +\sum_{i=1}^n
	{\Phi}^{\mu V}_{\lambda}{}_{\overline{i}} (z)
	\otimes v_{\overline{i}} , $$
   For the type II, the components are defined similarly.

   We take the following normalization.
\begin{equation}\label{E:norm}
\Phi^{\mu_j V}_{\mu_i}(z)|\mu_i>=|\mu_j>\otimes v_{ij}+ \mbox{terms of positive
powers of $z$},
\end{equation}
where $v_{12}=v_{\overline 1}, v_{21}=v_1, v_{34}=v_{n},
v_{43}=v_{\overline n}$.

     We also introduce the elements $a_{\overline{1}}(k)$ in the Heisenberg
algebra such that
     $$[a_j(k), a_{\overline{1}}(l)]=\delta_{j1}\delta_{k, -l}$$
where $$a_{\overline 1}(k)=\sum_{i=1}^n\frac{k}{[k/2]^2}
\frac{[(n+1-i)k][(n+1)k]_1}{[(n+1)k][(n+1-i)k]_i}a_i(k).$$

   \begin{theorem}
   1) For the type one vertex operators associated to $(\lambda, \mu)=
   (\mu_1, \mu_2)$, $(\mu_2, \mu_1)$, $(\mu_3, \mu_4)$, and $(\mu_4, \mu_3)$
   respectively, one has
   \begin{eqnarray} \nonumber
      &&\widetilde{\Phi}^{\mu_j V}_{\mu_i}{}_{\overline 1}(z)
      =exp(\sum_{k=1}^{\infty}\frac{[k/2]}kq^{(n+1/4-\delta_{n1})k}
a_{\overline 1}(-k)z^k)
      exp(\sum_{k=1}^{\infty}\frac{[k/2]}kq^{-(n+3/4-\delta_{n1})k}
a_{\overline 1}(k)z^{-k})\\
    &&\qquad\qquad
  e^{\lambda_1}(q^{(n+1/2)}z)^{-\lambda_1(0)+1-(\lambda_1|\mu_i)}\partial
Z^+_1(q^{n+1/2}z)c_{ij},
\label{E:vo}\\
      &&{\Phi}^{\mu V}_{\lambda}{}_{\overline {j+1}}(z)
      =-[{\Phi}^{\mu V}_{\lambda}{}_{\overline j}(z),
	f_j]_{q_j}, \quad
      {\Phi}^{\mu V}_{\lambda}{}_j(z)
      =[{\Phi}^{\mu V}_{\lambda}{}_{j+1}(z),
	f_j]_{q_j} \nonumber\\
      &&{\Phi}^{\mu V}_{\lambda}{}_{n}(z)
      =[{\Phi}^{\mu V}_{\lambda}{}_{\overline n}(z),
	f_j]_{q_j}, \nonumber
   \end{eqnarray}
where $c_{ij}$ are constants
for the four cases $(\mu_i, \mu_j)$ with $c_{12}=1$.

   2) The type two vertex operators are given by
   \begin{eqnarray*}
      {\Phi}^{V\mu}_{\lambda}{}_1(z)
     &=&exp(-\sum_{k=1}^{\infty}\frac{[k/2]}kq^{-k/4}
a_{\overline 1}(-k)z^k)
      exp(-\sum_{k=1}^{\infty}\frac{[k/2]}kq^{-3k/4}
a_{\overline 1}(k)z^{-k})c_{ij}'\\
     &&\qquad\qquad
e^{-\lambda_1}(q^{-1/2}z)^{\lambda_1(0)+1+(\lambda_1|\mu_i)}
Z^+_1(q^{-1/2}z)c_{ij}'
,\\
{\Phi}^{\mu V}_{\lambda}{}_{\overline {j}}(z)
      &=&-[{\Phi}^{\mu V}_{\lambda}{}_{\overline {j+1}}(z),
	e_j]_{q_j}, \quad
      {\Phi}^{\mu V}_{\lambda}{}_{j+1}(z)
      =[{\Phi}^{\mu V}_{\lambda}{}_{j}(z),
	e_j]_{q_j},\\
{\Phi}^{\mu V}_{\lambda}{}_{\overline n}(z)
     &=&[{\Phi}^{\mu V}_{\lambda}{}_{n}(z),
	e_n]_{q},
   \end{eqnarray*}
where $c'_{ij}$ are constants
for the four cases $(\mu_i, \mu_j)$ with $c'_{12}=1$.
  \end{theorem}
  \bigskip
  \noindent
  {\it Proof.}
  The proof of intertwining relations will be given
  in Section 4 and 5 on the space $\tilde{\cal F}$. We will only give the
argument for the type I case, and type II case can be treated similarly
word by word as in case I.
  It is obvious that each components of the vertex operators
  commute or anticommute with $Q^-_j$, thus the results are passed to the
  irreducible modules $V(\mu_i)$.
  \hspace{\fill} $\Box$

\section{Analysis of intertwining relations}

The vertex operator $\Phi(z)=\sum_{i=1}^n\Phi_i(z)\otimes v_i
+\sum_{i=1}^n\Phi_{\overline i}(z)\otimes v_i$ is determined by the
normalization (\ref{E:norm}) and the following relations: for $i=1, \cdots , n$
\begin{eqnarray}\label{E:R0}
\ t_j\Phi_i(z)t_j^{-1}&=& q_j^{\delta_{i,j+1}-\delta_{ij}}
\Phi_i(z), \quad j=0, 1, \cdots, n\\
t_j\Phi_{\overline i}(z)t_j^{-1}&=&q_j^{\delta_{ij}-\delta_{i,j+1}}
\Phi_{\overline i}(z), \quad j=0, 1, \cdots, n\label{E:R00}\\
\ [\Phi_i(z), e_j]&=&t_j\Phi_{j+1}(z)\delta_{ij}, \ \ j=1, \cdots, n-1
\label{E:R1}\\
\ [\Phi_{\overline i}(z), e_j]&=&-t_j\Phi_{\overline j}(z)\delta_{i,j+1},
\ \ j=1, \cdots, n\label{E:R2}\\
\ [\Phi_{\overline i}(z), e_0]&=&zt_0\Phi_{1}(z)\delta_{i1}, \ \
[\Phi_i(z), e_0]=0, \ \ j=1, \cdots n,
\label{E:R3}\\
\ [\Phi_{i}(z), e_n]&=&t_n\Phi_{\overline n}(z)\delta_{in}, \ \
[\Phi_{\overline i}(z), e_n]=0, \ \ j=1, \cdots n,
\label{E:R4}\\
\ [\Phi_i(z), f_j]_{q_j^{\delta_{i,j+1}-\delta_{ij}}}&=&\Phi_j(z)
\delta_{i,j+1},
\ \ i, j=1, \cdots n\label{E:R5}\\
\ [\Phi_{\overline i}(z), f_j]_{q_j^{\delta_{ij}-\delta_{i,j+1}}}&=
&-\Phi_{\overline{j+1}}(z)\delta_{ij},
\ \ i, j=1, \cdots n\label{E:R6}\\
\ [\Phi_i(z), f_0]_{q^{\delta_{i1}}}
&=&z^{-1}\Phi_{\overline 1}(z)\delta_{i1}, \quad
[\Phi_{\overline i}(z), f_0]_{q^{-\delta_{i1}}}=0, \label{E:R7}\\
\ [\Phi_{\overline i}(z), f_n]_{q^{\delta_{in}}}&=
&-\Phi_{n}(z)\delta_{in}, \quad [\Phi_{\overline i}(z), f_n]_{q^{\delta_{in}}}
=0
\label{E:R8}
\end{eqnarray}

With the construction of $\Phi_{\overline 1}(z)$ given in
(\ref{E:vo}), we define the other components:
\begin{eqnarray} \label{E:D1}
\Phi_{\overline{i+1}}(z)&=&-[\Phi_{\overline {i}}(z), f_i]_{q^{1/2}},
\ \ i=1, \cdots, n-1\\
\Phi_{n}(z)&=&[\Phi_{\overline {n}}(z), f_n]_{q}, \label{E:D2}\\
\Phi_{i}(z)&=&[\Phi_{i+1}(z), f_i]_{q^{1/2}},
\ \ i=1, \cdots, n-1 \label{E:D3}
\end{eqnarray}

In order to study $q$-commutators, we recall the following identities from
\cite{kn:J}:
\begin{eqnarray} \label{E:Adj1}
\ \left[a, [b, c]_u\right]_v&=&\left[[a,b]_x,c\right]_{uv/x}+x\,
\left[b, [a,c]_{v/x}\right]_{v/x},
\qquad x\neq 0\\
\ \left[[a, b]_u, c\right]_v&=&\left[a,[b, c]_x\right]_{uv/x}+x\,
\left[[a, c]_{v/x}, b\right]_{u/x},
\qquad x\neq 0  \label{E:Adj2}
\end{eqnarray}

\begin{prop} Let $\Phi_i(z)$ be vertex operators defined by
$\Phi_{\overline 1}(z)$ via (\ref{E:D1}, \ref{E:D2}, \ref{E:D3})
and satisfy the relations:
\begin{eqnarray*}
\ [\Phi_{\overline 1}(z), f_j]_{q_j^{\delta_{1j}}}&=&-\Phi_{\overline 2}(z)
\delta_{1j}, \quad
[\Phi_{\overline 1}(z), e_j]=0, \quad j=1, \cdots n\\
\ \left[f_1, [f_1, \Phi_{\overline 1}(z)]_{q^{-1/2}}\right]_{q^{1/2}}&=&0.
\end{eqnarray*}
Then we have for $i=0, 1, \cdots, n$, $j=1, \cdots, n$
\begin{eqnarray}
\ [\Phi_i(z), e_j]&=&t_j\Phi_{j+1}(z)\delta_{ij}, \label{E:R10}\\
\ [\Phi_{\overline i}(z), e_j]&=&-t_j\Phi_{\overline j}(z)\delta_{i,j+1},
\label{E:R11}\\
\ [\Phi_{i}(z), e_n]&=&t_n\Phi_{\overline n}(z)\delta_{in}, \ \
[\Phi_{\overline i}(z), e_n]=0, \label{E:R12}\\
\ [\Phi_i(z), f_j]_{q_j^{\delta_{i,j+1}-\delta_{ij}}}&=&
\Phi_j(z)\delta_{i,j+1}, \label{E:R13}\\
\ [\Phi_{\overline i}(z), f_j]_{q_j^{\delta_{ij}-\delta_{i,j+1}}}&=
&-\Phi_{\overline{j+1}}(z)\delta_{i,j+1}, \label{E:R14}
\end{eqnarray}
where  we identify $\Phi_{n+1}(z)=\Phi_{\overline n}(z)$ in (\ref{E:R10}),
$\Phi_{\overline 0}(z)=-z\Phi_{1}(z)$ in (\ref{E:R11}),
$\Phi_{\overline{n+1}}(z)=\Phi_{n}(z)$ in (\ref{E:R13}).
\end{prop}
{\it Proof}. Using relations of the quantum affine algebras $U_q(C_n^{(1)})$,
we have for $j\neq i-1$:
\begin{eqnarray*}
\ [\Phi_{\overline{i}}(z), e_j]&=&-\left[[\Phi_{\overline{i-1}}(z), f_{i-1}]_
{q^{1/2}}, e_j\right]\\
\ &=&-\left[ [\Phi_{\overline{i-1}}(z), e_j], f_{i-1}\right]_{q^{1/2}}=0 \qquad
\mbox{by induction}
\end{eqnarray*}
\begin{eqnarray*}
\ [\Phi_{\overline{i}}(z), e_{i-1}]
&=&-\left[[\Phi_{\overline{i-1}}(z), f_{i-1}]_
{q^{1/2}}, e_{i-1}\right]\\
\ &=&\left[ \Phi_{\overline{i-1}}(z),
\frac{t_{i-1}-t_{i-1}^{-1}}{q^{1/2}-q^{-1/2}}
\right]_{q^{1/2}}=-t_{i-1}\Phi_{i-1}(z).
\end{eqnarray*}
\begin{eqnarray*}
\ [\Phi_{i}(z), e_j]&=&\left[[\Phi_{i+1}(z), f_{i}]_
{q^{1/2}}, e_j\right]=-\left[ \Phi_{\overline{i+1}}(z),
\frac{t_{i}-t_{i}^{-1}}{q_i-q_i^{-1}}\right]_{q_i}\delta_{ij}\\
&=&t_{j}\Phi_{j+1}(z)\delta_{ij},
\end{eqnarray*}
where we used induction to get $[\Phi_{i+1}, e_j]=0$ for $j\neq i$ in
the second equality.
Similarly we can prove the commutation relations of $\Phi_i(z)$ and $e_n$
for $i=n, \overline n$.

Now we consider the relations between $\Phi_i(z)$ and $f_j$. It follows from
(\ref{E:Adj2}) that
\begin{eqnarray*}
\ [\Phi_i(z), f_j]&=&\left[[\Phi_{i-1}(z), f_{i-1}]_{q^{1/2}}, f_j\right]\\
&=&\left[[\Phi_{i-1}(z), f_j], f_{i-1}\right]_{q^{1/2}}+
\left[\Phi_{i-1}(z), [f_{i-1}, f_j]\right]_{q^{1/2}}
\end{eqnarray*}
Then the relation (\ref{E:R5}) is true: when $i=j$ is the definition;
$i=j+1$ follows from Proposition \ref{prop2}, and other cases are
derived from the relation.

\begin{eqnarray*}
\ [\Phi_{\overline i}(z), f_{i-1}]_{q^{-1/2}}
&=&-\left[[\Phi_{\overline{i-1}}(z), f_{i-1}]_{q^{1/2}},
f_{i-1}\right]_{q^{-1/2}}\\
&=&-\left[\left[[\Phi_{\overline{i-2}}(z), f_{i-2}]_{q^{1/2}},
f_{i-1}\right]_{q^{1/2}}, f_{i-1}\right]_{q^{-1/2}}
\end{eqnarray*}
Now let use induction on $i$. The initial step $i=2$ is Proposition
(\ref{prop2}). From the induction hypothesis
$[\Phi_{\overline{i-2}}(z), f_{i-1}]_{q^{1/2}}=0$, we see the above
commutation turns to:
\begin{eqnarray*}
\ [\Phi_{\overline i}(z), f_{\overline{i-1}}]_{q^{-1/2}}&=&
-[\Phi_{\overline{i-2}}(z), \left[[f_{i-2}, f_{i-1}]_{q^{1/2}},
f_{i-1}\right]_{q^{-1/2}}]_{q^{1/2}}\\
&=&0. \qquad\mbox{by Serre relation}
\end{eqnarray*}
Clearly we have
$$
[\Phi_{\overline i}(z), f_{i+1}]=
-\left[[\Phi_{\overline{i-1}}(z), f_{i-1}]_{q^{1/2}}, f_{i+1}\right]=0.
$$
\begin{eqnarray*}
\ [\Phi_{\overline i}(z), e_0]&=&-[\Phi_{\overline{i-1}}(z), e_0], f_{i-1}]\\
&=&-\delta_{i2}[t_0\Phi_1(z), f_1]_{q^{1/2}}
=-\delta_{i2}t_0[\Phi_1, f_1]_{q^{-1/2}}\\
&=& 0  \qquad\mbox{by (\ref{E:R5})}
\end{eqnarray*}
where we used induction on $i$ in the second equality. Similarly we can
check that \linebreak
$[\Phi_j(z), f_0]=0$
for $j\neq 1$ by induction.

Next we claim that $[\Phi_{\overline 1}(z), e_0]=zt_0\Phi_1(z)$
is equivalent to $[\Phi_1(z), f_0]_{q^{-1}}=\Phi_{\overline 1}(z)z^{-1}$.
In fact, given the former we have that
\begin{eqnarray*}
\ [\Phi_1(z), f_0]&=&z^{-1}\left[
t_0^{-1}[\Phi_{\overline 1}(z), e_0], f_0\right]\\
&=&z^{-1}t_0^{-1}\left[[\Phi_{\overline 1}(z), e_0], f_0\right]_{q^{-1}}\\
&=&z^{-1}t_0^{-1}\left[\Phi_{\overline 1}(z), \frac{t_0-t_0^{-1}}{q-q^{-1}}
\right]_{q^{-1}}=z^{-1}\Phi_{\overline 1}(z).
\end{eqnarray*}

With the relations of
$[\Phi_{\overline i}(z), f_j]_{q_j^{\delta_{ij}-\delta_{i, j+1}}}=
-\Phi_{\overline{j+1}}(z)\delta_{ij}$ in hand, we can compute that
\begin{eqnarray*}
\ [\Phi_i(z), f_{i}]_{q^{-1/2}}&=&
\left[[\Phi_{i+2}(z), f_{i+1}]_{q_{i+1}}, f_i]_{q^{1/2}}, f_i
\right]_{q^{-1/2}} \\
&=&\left[\Phi_{i+2}(z), \left[[f_{i+1}, f_i]_{q^{1/2}}, f_i\right]_{q^{-1/2}}
\right]_{q_{1/2}}\\
&=&0
\end{eqnarray*}
where we have used induction on $n-i$ and $i\leq n-2$. When $i=n, n-1$
we can compute similarly. For example,
\begin{eqnarray*}
\ [\Phi_n(z), f_{n}]_{q^{-1}}&=&
\left[[\Phi_{\overline{n-1}}(z), f_{n-1}]_{q^{1/2}}, f_n]_{q}, f_n
\right]_{q^{-1}} \\
&=&\left[\Phi_{\overline{n-1}}(z), \left[[f_{n-1}, f_n]_{q}, f_n
\right]_{q^{-1}}
\right]_{q_{1/2}}\\
&=&0
\end{eqnarray*}

The other cases of $[\Phi_{i}(z), f_j]_{q_j^{-\delta_{ij}+\delta_{i,j+1}}}=
\Phi_{j}(z)\delta_{i-1,j}$ are shown similarly.

\begin{prop} \label{prop1}
The operator $\Phi(z)$ (cf. (\ref{E:vo})) satisfy the following relations:
\begin{eqnarray*}
\ [\Phi_{\overline 1}(z), X_i^+(w)]&=& 0, \ \ i=1, \cdots n\\
\ [\Phi_{\overline 1}(z), X_i^-(w)]&=& 0, \ \ i=2, \cdots n\\
\ [\Phi_{\overline 1}(z), X_1^-(w)]_{q^{-1/2}}&=&\frac zw q^n
  [\Phi_{\overline 1}(z), X_1^-(w)]_{q^{1/2}}
\end{eqnarray*}
\end{prop}
{\it Proof}. The first set of relations follow from our construction.
We only need to see $i=1$ in the second set. Writing
$\Phi_{\overline 1}(z)=\frac 1{(q^{1/2}-q^{-1/2})z}(\Phi_{\overline 1+}(z)
-\Phi_{\overline 1-}(z))$, we have
\begin{eqnarray}
\Phi_{\overline 1\ep}(z)X_{1\ep'}^+(z)&=&:\Phi_{\overline 1\ep}(z)
X_{1\ep'}^+(z):
q^{\ep/2}\frac{q^{n+1/2}z-q^{\ep'/2-\ep/2}w}{q^{n+1/2}z-w}
\label{E:contr1}\\
X_{1\ep'}^-(z)\Phi_{\overline 1\ep}(z)&=&:\Phi_{\overline
1\ep}(z)X_{1\ep'}^-(z):
q^{\ep/2}\frac{q^{\ep'/2-\ep/2}w-q^{n+1/2}z}{q^{n+1/2}w-z},
\label{E:contr2}
\end{eqnarray}

\begin{eqnarray*}
&&[\Phi_{\overline 1}(z), X_i^+(w)] \\
&=&\sum_{\ep=\pm 1} [\Phi_{\overline 1\ep}(z), X_{i,-\ep}^+(w)]\\
&=&\sum_{\ep=\pm 1}:\Phi_{\overline 1\ep}(z)X_{1\ep'}^-(z):
\frac{q^{\ep/2}(q^{n+1/2}z-q^{-\ep}w)}w \delta(\frac{q^{n+1/2}z}w)\\
&=&0
\end{eqnarray*}
since $:\Phi_{\overline 1+}(z)X_{1-}^+(w):=
:\Phi_{\overline 1-}(z)X_{1-}^+(w):$.

The third relation is another form of the following identity.
$$
\begin{array}{ll}
&q^{1/2}(w-q^nz)\Phi_{\overline 1}(z)X_1^-(w)+
(q^{n+1}z-w)X_1^-(w)\Phi_{\overline 1}(z)\\
&=\sum_{\ep', \ep=\pm}
q^{1/2}(w-q^nz)\Phi_{\overline 1\ep}(z)X_{1\ep'}^-(w)+
(q^{n+1}z-w)X_{1\ep'}^-(w)\Phi_{\overline 1\ep}(z)\\
&=\sum_{\ep', \ep=\pm}:\Phi_{\overline 1\ep}(z)X_{1\ep'}^-(w):
q^{-\ep/2}(q^{n+1}z-w)(w-q^nz)\delta(\frac{q^{(2n+\ep+1)/2}z}w)
q^{-n-\ep/2}z^{-1}\\
&=0
\end{array}
$$
\hspace{\fill} $\Box$

In particular, we have
\begin{equation}\label{E:1to0}
[\Phi_{\overline 1}(z), X_1^-(1)]_{q^{-1/2}}
=zq^n[\Phi_{\overline 1}(z), X_1^-(0)]_{q^{1/2}}
\end{equation}

\begin{prop}\label{prop2}
The operator constructed in (\ref{E:vo}) satisfies the following
Serre-like relation:
$$\Phi_{\overline 1}(z)f_1^2-(q^{1/2}+q^{-1/2})f_1\Phi_{\overline 1}(z)f_1
+f_1^2\Phi_{\overline 1}(z)=0$$
\end{prop}
{\it Proof}. It follows from (\ref{E:contr1}, \ref{E:contr2}) that
\begin{eqnarray*}
\Phi_{\overline 1\ep}(w)X_{1\ep_1}^-(z_1)X_{1\ep_2}^-(z_2)&=&
:\Phi_{\overline 1\ep}(w)X_{1\ep_1}^-(z_1)X_{1\ep_2}^-(z_2):\\
&&\prod_{i=1,2}\frac{q^{n+1/2}w-q^{-1/2}z_i}{q^{(2n+\ep+1)/2}w-z_i}
\cdot \frac{q^{\ep_1/2}z_1-q^{\ep_2/2}z_2}{z_1-q^{-1}z_2}\\
 X_{1\ep_1}^-(z_1)\Phi_{\overline 1\ep}(w)X_{1\ep_2}^-(z_2)&=&
:\Phi_{\overline 1\ep}(w)X_{1\ep_1}^-(z_1)X_{1\ep_2}^-(z_2):\\
&&\frac{z_1-q^nw}{z_1-q^{(2n+\ep+1)/2}w}\dot
\frac{q^{n+1/2}w-q^{-1/2}z_2}{q^{(2n+\ep+1)/2}w-z_2}
\cdot \frac{q^{\ep_1/2}z_1-q^{\ep_2/2}z_2}{z_1-q^{-1}z_2}\\
X_{1\ep_1}^-(z_1)X_{1\ep_2}^-(z_2)\Phi_{\overline 1\ep}(w)&=&
:\Phi_{\overline 1\ep}(w)X_{1\ep_1}^-(z_1)X_{1\ep_2}^-(z_2):\\
&&\prod_{i=1,2}\frac{z_i-q^nw}{z_i-q^{(2n+\ep+1)/2}w}
\cdot \frac{q^{\ep_1/2}z_1-q^{\ep_2/2}z_2}{z_1-q^{-1}z_2}
\end{eqnarray*}
Then we have
\begin{eqnarray*}
&&\Phi_{\overline 1\ep}(w)X^-_{1\ep_1}(z_1)X^-_{1\ep_2}(z_2)
-(q^{1/2}+q^{-1/2})X^-_{1\ep_1}(z_1)\Phi_{\overline 1}(w)X_{1\ep_2}(z_2)
+X^-_{1\ep_1}(z_1)X^-_{1\ep_2}(z_2)\Phi_{\overline 1}(z)\\
&=&:\Phi_{\overline 1\ep}(w)X^-_{1\ep_1}(z_1)X^-_{1\ep_2}(z_2):
\frac{q^{\ep_1/2}z_1-q^{\ep_2/2}z_2}{(z_1-q^{-1}z_2)
\prod_{i=1,2}(q^{(2n+\ep+1)/2}w-z_i)}\\
&&\  \cdot \left\{ (q^{n+1/2}w-q^{-1/2}z_1)
(q^{n+1/2}w-q^{-1/2}z_2)-(q^{1/2}+q^{-1/2})\cdot\right.\\
&&\ \ \ \left. \cdot(z_1-q^nw)(q^{-1/2}z_2-q^{n+1/2}w)+(z_1-q^nw)(z_2-q^nw)
\right\}\\
&=&:\ \ :
\frac{q^{\ep_1/2}z_1-q^{\ep_2/2}z_2}{
\prod_{i=1,2}(q^{(2n+\ep+1)/2}w-z_i)}
q^n(1-q)w
\end{eqnarray*}
Note that the contraction function is antisymmetric respect to $z_1
\mapsto z_2$ when $\ep_1=\ep_2$. Therefore,
\begin{eqnarray*}
\lefteqn{Sym_{z_1, z_2}\Phi_{\overline 1\ep}(w)X^-_{1\ep_1}(z_1)
X^-_{1\ep_1}(z_2)}\\
&& -(q^{1/2}+q^{-1/2})X^-_{1\ep_1}(z_1)\Phi_{\overline 1}(w)X^-_{1\ep_1}(z_2)
+X^-_{1\ep_1}(z_1)X^-_{1\ep_1}(z_2)\Phi_{\overline 1}(z)=0
\end{eqnarray*}
Furthermore when $\ep_1=-\ep_2$ we have
\begin{eqnarray*}
&&\sum_{\ep_1}\Phi_{\overline 1\ep}(w)X^-_{1\ep_1}(z_1)X^-_{1-\ep_1}(z_2)
-(q^{1/2}+q^{-1/2})X^-_{1\ep_1}(z_1)\Phi_{\overline 1}(w)X_{1-\ep_1}(z_2)+\\
&&\qquad\qquad X^-_{1\ep_1}(z_1)X^-_{1-\ep_1}(z_2)\Phi_{\overline 1}(z)\\
&=&:\Phi_{\overline 1\ep}(w)X^-_{1\ep_1}(z_1)X^-_{1-\ep_1}(z_2):
\frac{(q^{1/2}+q^{-1/2})(z_1-z_2)}{
\prod_{i=1,2}(q^{(2n+\ep+1)/2}w-z_i)}
q^n(1-q)w
\end{eqnarray*}
which is also antisymmetric with respect to $z_1, z_2$. We thus have that
\begin{eqnarray*}
\lefteqn{Sym_{z_1, z_2}\Phi_{\overline 1}(w)X^-_{1}(z_1)
X^-_{1}(z_2)}\\
&& -(q^{1/2}+q^{-1/2})X^-_{1}(z_1)\Phi_{\overline 1}(w)X_{1}(z_2)
+X^-_{1}(z_1)X^-_{1}(z_2)\Phi_{\overline 1}(z)=0
\end{eqnarray*}
Picking up coefficients of $(z_1z_2)^m$ we proved the identity.
\hspace{\fill} $\Box$

\section{Proof of $[\Phi_{\overline 1}(z), e_0]=zt_0\Phi_{1}(z)$}

In the last section we have shown that we are left to check only the
relation in the section title. The key of the the proof
 are the $q$-commutation
relations (\ref{E:Adj1}) and (\ref{E:Adj2}) (cf. \cite{kn:J}).

As in \cite{kn:J} we introduce the
 twisted commutators $[b_1, \cdots, b_n]_{v_1\cdots v_{n-1}}$ and
$[b_1, \cdots, b_n]'_{v_1\cdots v_{n-1}}$
inductively by $[b_1, b_2]_v=[b_1, b_2]'_v=b_1b_2-vb_2b_1$ and
\begin{eqnarray*}
\ [b_1, \cdots, b_n]_{v_1\cdots v_{n-1}}&=&\left[b_1,
[b_2, \cdots, b_n]_{v_1\cdots v_{n-2}}\right]_{v_{n-1}}\\
\ [b_1, \cdots, b_n]'_{v_1\cdots v_{n-1}}&=&\left[[b_1,
\cdots, b_{n-1}]'_{v_1\cdots v_{n-2}}, b_n\right]_{v_{n-1}}
\end{eqnarray*}

In this notation we have that
\begin{eqnarray*}
&&e_0\\
&=&[X_1^-(0), \cdots, X_n^-(0),
 X_{n-1}^-(0), \cdots, X_2^-(0),
X_1^-(1)]_{q^{-1/2}\cdots q^{-1/2}\, q^{-1}\, q^{-1/2}\cdots q^{-1/2}\, 1}
\gamma K_{\theta}^{-1}\frac{(-1)^n}{[2]_1}.
\end{eqnarray*}

Write $\hat{e_0}= [2]_1e_0\gamma^{-1}K_{\theta}$, then the relation we want
is
$$[\Phi_{\overline 1}(z), \hat{e_0}]_{q^{-1}}=(-1)^n(1+q)z\Phi_1(z). $$

\begin{lemma}\label{L10}
\begin{eqnarray*}
\ &&[X_n^-(0), X_{n-1}^-(0), \cdots, X_2^-(0), X_1^-(1)]_{q^{-1/2}\,
\cdots q^{-1/2}\, q^{-1}}\\
&=&\ \ q^{-n/2}[X_n^-(1), X_{n-1}^-(0), \cdots, X_2^-(0), X_1^-(0)]_{q^{1/2}\,
\cdots q^{1/2}\, q}\\
&=&\ \ (-1)^{n-1}[X_1^-(0), X_{2}^-(0), \cdots, X_{n-1}^-(0), X_n^-(1)
]_{q^{-1}\, q^{-1/2}
\cdots q^{-1/2}}
\end{eqnarray*}
\end{lemma}
{\it Proof}. From the Serre relations
$[X_i^{-}(0), X_j^-(1)]_{q^{(\a_i, \a_j)}}=
-[X_j^{-}(0), X_i^-(1)]_{q^{(\a_i, \a_j)}}$, it follows that
\begin{eqnarray*}
\ &&[X_n^-(0), X_{n-1}^-(0), \cdots, X_2^-(0), X_1^-(1)]_{q^{-1/2}\,
\cdots q^{-1/2}\, q^{-1}}\\
&=&-[X_n^-(0), X_{n-1}^-(0), \cdots, X_3^-(0), X_1^-(0), X_2^-(1)]_{q^{-1/2}\,
\cdots q^{-1/2}\, q^{-1}}\\
&=&(-1)^2[X_n^-(0), X_{n-1}^-(0), \cdots, X_{1}^-(0), X_2^-(0), X_3^-(1)]
]_{q^{-1/2}
\cdots q^{-1/2}\, q^{-1}} \quad\mbox{inductively}\\
&=&(-1)^{n-1}[X_1^-(0), X_{2}^-(0), \cdots, X_{n-1}^-(0), X_n^-(1)
]_{q^{-1}\, q^{-1/2}
\cdots q^{-1/2}}\\
&=& q^{-1/2}(-1)^{n-2}
[X_1^-(0), X_{2}^-(0), \cdots, X_{n}^-(1), X_{n-1}^-(0)
]_{q\, q^{-1/2}
\cdots q^{-1/2}}\\
&=& q^{-n/2}[X_n^-(1), X_{n-1}^-(0), \cdots, X_2^-(0), X_1^-(0)]_{q^{1/2}
\cdots q^{1/2}\, q}
\end{eqnarray*}

\begin{lemma} \label{L11} For $i\leq n-1$ we have
\begin{eqnarray*}
&&\left[[X_1^-(0), X^-_2(0), \cdots, X_i^-(0)]'_{q^{-1/2}\cdots q^{-1/2}},
\right. \\
&& \quad \left. [X_{i+1}^-(0), X_i^-(0), \cdots, X_2^-(0),
X_1^-(1)]_{q^{-1/2}\cdots q^{-1/2}}\right]_{q^{1/2}}=0\\
\end{eqnarray*}
\end{lemma}
{\it Proof}. Following \cite{kn:Dr},
for $y\in End(\tilde{\cal F})$ we define $q$-adjoint operators
$ad\, X_i^{\pm}(0)=S_i^{\pm}$ by
$$S_i^{\pm}(y)=X_i^{\pm}(0)y-K_i^{\pm 1}yK_i^{\mp 1}X_i^{\pm}(0).$$
We can directly check by Serre relations that
\begin{eqnarray}
&[S_i^{\pm}, S_j^{\pm}]=0, &\mbox{if $A_{ij}=0$} \label{E:Adj3}\\
&[S_i^{\pm}, S_i^{\pm}, S_j^{\pm}]_{q_i, q_i^{-1}}=0, & \mbox{if $A_{ij}=-1$}
\label{E:Adj4}\\
& S_i^{\pm 2}(X_j^{\pm}(m))=0 &\mbox{if $A_{ij}=-1$} \label{E:Adj5}
\end{eqnarray}
We now prove the identity by induction on $i$. $i=1$ is the Serre
relation:
$$[X_1^-(0), X_2^-(0), X_1^-(1)]_{q^{-1/2}q^{1/2}}
= -[X_1^-(0), X_1^-(0), X_2^-(1)]_{q^{-1/2}q^{1/2}}=0.
$$
Set
\begin{eqnarray*}
&A'=&[X_1^-(0), X^-_2(0), \cdots, X_{i-1}^-(0)
]'_{q^{-1/2}\cdots q^{-1/2}}\\
&A=&[X_{i-1}^-(0), \cdots, X_2^-(0),
X_1^-(1)]_{q^{-1/2}\cdots q^{-1/2}}=S_{i-1}^-\cdots S_2^-(X_1^-(1)).
\end{eqnarray*}
The inductive assumption is then
$$
\left[ A', [X_i^-(0), A]_{q^{-1/2}}\right]_{q^{1/2}}=0.
$$
Then we have
\begin{eqnarray}\nonumber
\ &&[X_i^-(0), X_{i+1}^-(0), X_i^-(0), A]_{q^{-1/2}\, q^{-1/2}\, 1}=
S_i^-S_{i+1}S_i^-S_{i-1}^-\cdots S_2^-(X_2^-(1))\\
  &=&\frac 1{[2]_1}\left(S_i^{-2}S_{i+1}^-S_{i-1}^-\cdots S_2^-(X_2^-(1))
  +S_{i+1}^{-}(S_{i}^{-2}S_{i-1}^-)S_{i-2}^-\cdots S_2^-(X_2^-(1))\right)
  \nonumber\\
  &=&\frac 1{[2]_1}S_{i+1}^-\left([2]_1S_i^-S_{i-1}^-S_i^-S_{i-2}^-\cdots
  S_2^-(X_2(1))-S_{i-1}^-S_{i-1}^{-}S_i^{-2}S_{i-2}^-\cdots S_2^-(X_1^-(1))
  \right)\nonumber\\
  &=& 0 \label{E:Adj6}
\end{eqnarray}
where we have used (\ref{E:Adj3}-\ref{E:Adj4}). From (\ref{E:Adj2}) and $
[A', X_{i+1}^-(0)]=0$ it
follows that
\begin{eqnarray*}
 &&\left[[A', X_i^-(0)]_{q^{-1/2}}, [X_{i+1}^-(0), X_i^-(0),
A]_{q^{-1/2}\, q^{-1/2}}\right]_{q^{1/2}}\\
&=&\left[A', [X_i^-(0), X_{i+1}^-(0), X_i^-(0), A]_{q^{-1/2}\, q^{-1/2}\, 1}
\right]+\\
&&\left[[X_{i+1}^-(0),
A', X_i^-(0), A]_{q^{-1/2}\, q^{1/2}\, q^{-1/2}},
X_i^-(0)\right]_{q^{-1/2}}\\
&=& 0\qquad\mbox{by (\ref{E:Adj6}) and induction hypothesis}
\end{eqnarray*}
\hspace{\fill} $\Box$

\begin{lemma}\label{L11a} On the space $End(\tilde{\cal F})$ we have
\begin{equation}\nonumber
\left[[\Phi_{\overline 1}(z), X_1^-(0), \cdots X_{n-1}^-(0)]'_{q^{-1/2}\cdots
q^{-1/2}},
[X_{n-1}^-(0), \cdots, X_1^-(0)]_{q^{1/2}\cdots q^{1/2}}
\right]_{q^{1/2}}=0
\end{equation}
\end{lemma}
{\it Proof}. The proof is by induction and similar to that of
Lemma \ref{L11}, but with the $q$-adjoint operators $T_i^{\pm}$:
$T_i^{\pm}(y)=
X^{\pm}_i(0)y-K_i^{\mp 1}yK_i^{\pm 1}X_i^{\pm}(0)$.
The operators $T_i^{\pm}$ satisfy the same relations (\ref{E:Adj3},
\ref{E:Adj4}, \ref{E:Adj5}) as the operators $S_i^{\pm}$.
Then we can express the commutators in the identities as follows:
\begin{eqnarray*}
&[X_{n-1}^-(0), \cdots, X_1^-(0), \Phi_{\overline 1}(z)]_{q^{1/2}\cdots
q^{1/2}} &=T_{n-1}^-\cdots T_1^-(\Phi_{\overline 1}(z)), \\
&[X_{n-1}^-(0), \cdots, X_2^-(0), X_1^-(0)]_{q^{1/2}\cdots q^{1/2}}
 &=T_{n-1}^-\cdots T_2^-(X_1^-(0)).
\end{eqnarray*}

We claim that for $i\leq n-1$
\begin{eqnarray}\nonumber
&&T_i^-T_{i+1}^-T_i^-\cdots T_1^-(\Phi_{\overline 1}(z))\\
&=& \left[ X_i^-(0), [X_{i+1}^-(0), \cdots, X_1^-(0),
\Phi_{\overline 1}(z)]_{q^{1/2}\cdots q^{1/2}}\right]_{q^{-\frac 12
\delta_{n-1,i}}} \nonumber\\
&=&0. \label{E:Adj7}
\end{eqnarray}
The case of $i=1$ is exactly the Serre-like relation (cf. Prop. \ref{prop2}):
$$
[X_1^-(0), X_1^-(0), \Phi_{\overline 1}(z)]_{q^{1/2}\, q^{-1/2}}=0.
$$
Noting that $T_i(\Phi_{\overline 1}(z))=0$ for $i\geq 2$, we have
\begin{eqnarray*}
&&T_i^-T_{i+1}^-T_i^-\cdots T_1^-(\Phi_{\overline 1}(z))\\
&&=\frac 1{[2]_1}(T_i^{-2}T_{i+1}^-+T_{i+1}^-T_i^{-2})T_{i-1}\cdots T_1^-
(\Phi_{\overline 1}(z))\\
&&=\frac 1{[2]_1}T_{i+1}^-([2]_1T_{i}^-T_{i-1}^-T_{i}^--T_{i-1}T_i^{-2})
T_{i-2}^-\cdots T_1^-(\Phi_{\overline 1}(z))\\
&&=0.
\end{eqnarray*}
where we have used induction on $i$. Then we have for
$B=[X_{n-2}^-(0), \cdots,
X_1^-(0)]_{q^{1/2}\cdots q^{1/2}}$
\begin{eqnarray*}
&& \left[[X_{n-1}^-(0), \cdots, X_1^-(0)]_{q^{1/2}\cdots q^{1/2}},
[X_{n-1}^-(0), \cdots, X_{1}^-(0), \Phi_{\overline 1}(z)]_{q^{1/2}\cdots
q^{1/2}}\right]_{q^{-1/2}}\\
&&=\left[X_{n-1}^-(0), \left[B,
[X_{n-1}^-(0), \cdots, X_{1}^-(0), \Phi_{\overline 1}(z)]_{q^{1/2}\cdots
q^{1/2}}\right]\right]+\\
&&\qquad \left.\left[X_{n-1}^-(0),
[X_{n-1}^-(0), \cdots, X_{1}^-(0), \Phi_{\overline 1}(z)]_{q^{1/2}\cdots
q^{1/2}}\right]_{q^{-1/2}}, B\right]_{q^{1/2}}\\
&&=0   \qquad\qquad\mbox{by (\ref{E:Adj2}) and (\ref{E:Adj7})}
\end{eqnarray*}
\hspace{\fill} $\Box$

\begin{lemma} \label{L12} For $i\leq n-2$ we have
\begin{eqnarray*}
&&\left[[\Phi_{\overline 1}(z), X_1^-(0), \cdots, X_i^-(0)]'_{q^{-1/2}\cdots
q^{-1/2}}, [X_{i+1}^-(0), \cdots, X_2^-(0), X_1^-(1)]_{q^{-1/2}\cdots
q^{-1/2}\, q^{-1/2}}\right]\\
&=& -zq^{n+(1-i)/2}\left[\Phi_{\overline{i+2}}(z),
[X_1^-(0), \cdots, X_i^-(0)]'_{q^{-1/2}\cdots q^{-1/2}}\right]_{q^{-1}}.
\end{eqnarray*}
\end{lemma}
{\it Proof}. From Lemma \ref{L11} and (\ref{E:R5}) it follows that
\begin{eqnarray*}
&&\left[[\Phi_{\overline 1}(z), X_1^-(0), \cdots, X_i^-(0)]'_{q^{-1/2}\cdots
q^{-1/2}}, [X_{i+1}^-(0), \cdots, X_2^-(0), X_1^-(1)]_{q^{-1/2}\cdots
q^{-1/2}\, q^{-1/2}}\right]\\
&=&\left[\left[\Phi_{\overline 1}(z), [X_1^-(0), \cdots,
X_i^-(0)]'_{q^{-1/2}\cdots
q^{-1/2}}\right]_{q^{-1/2}}, [X_{i+1}^-(0), \cdots, X_2^-(0),
X_1^-(1)]_{q^{-1/2}\cdots q^{-1/2}}\right]\\
&=&\left[\Phi_{\overline 1}(z),
\left[[X_1^-(0), \cdots, X_i^-(0)]'_{q^{-1/2}\cdots q^{-1/2}},
[X_{i+1}^-(0), \cdots, X_2^-(0),
X_1^-(1)]_{q^{-1/2}\cdots q^{-1/2}}\right]_{q^{1/2}}
\right]_{q^{-1}}\\
&& +q^{1/2}\left[[\Phi_{\overline 1}(z), X_{i+1}^-(0), \cdots, X_2^-(0),
X_1^-(1)]_{q^{-1/2}\cdots q^{-1/2}}, [X_1^-(0), \cdots,
X_i^-(0)]'_{q^{-1/2}\cdots q^{-1/2}}\right]_{q^{-1}}\\
&=& q^{1/2}\left[[X_{i+1}^-(0), \cdots, X_2^-(0), \Phi_{\overline 1}(z),
X_1^-(1)]_{q^{-1/2}\cdots q^{-1/2}}, [X_1^-(0), \cdots,
X_i^-(0)]'_{q^{-1/2}\cdots q^{-1/2}}\right]_{q^{-1}}\\
&=&-q^{n+1/2}z\left[[X_{i+1}^-(0), \cdots, X_2^-(0), \Phi_{\overline 2}(z)
]_{q^{-1/2}\cdots q^{-1/2}}, [X_1^-(0), \cdots,
X_i^-(0)]'_{q^{-1/2}\cdots q^{-1/2}}\right]_{q^{-1}}
\end{eqnarray*}
where we used Lemma \ref{L11}, the fact of $\Phi_{\overline 1}(z)$
commuting with $X_i^-(0)$ for $i\geq 2$ (cf. \ref{E:R5}) and (\ref{E:1to0}).
Invoking the definition of $\Phi_{\overline i}(z)$ and converting the
$q$-commutators, we obtained the required identity.
\hspace{\fill} $\Box$

\begin{lemma}\label{L13} For $i\leq n-2$ we have
\begin{eqnarray*}
&&\left[X_i^-(0), \cdots X_n^-(0), \cdots, X_{i+1}^-(0),
[\Phi_{\overline 1}(z), X_1^-(0), \cdots,
X_{i-1}^-(0)]_{q^{-1/2}\cdots q^{-1/2}}, \right. \\
&&\qquad \left. [X_i^-(0), \cdots, X_2^-(0), X_1^-(1)]_{q^{-1/2}\cdots
q^{-1/2}}
\right]_{1\, q^{-1/2}\cdots q^{-1}\cdots q^{-1/2}\, 1}\\
&=&(-1)^{n-i-1}zq^{i/2+1}\left[X_i^-(0), \Phi_{i+1}(z), [X_1^-(0), \cdots
X_{i-1}^-(0)]'_{q^{-1/2}\cdots q^{-1/2}}\right]_{q^{-1}\, 1}\\
&=&(-1)^{n-i}zq^{(i+1)/2}
\left(\left[\Phi_i(z), [X_1^-(0), \cdots
X_{i-1}^-(0)]'_{q^{-1/2}\cdots q^{-1/2}}\right]_{q^{-1/2}}+\right.\\
&&\qquad\qquad \left. q^{1/2}\left[\Phi_{i+1}(z), [X_1^-(0), \cdots
X_{i}^-(0)]'_{q^{-1/2}\cdots q^{-1/2}}\right]_{q^{-1/2}}\right).
\end{eqnarray*}
\end{lemma}
{\it Proof}. It follows from Lemma \ref{L12} that
\begin{eqnarray*}
&&\left[X_i^-(0), \cdots X_n^-(0), \cdots X_{i+1}^-(0),
[\Phi_{\overline 1}(z), X_1^-(0), \cdots,
X_{i-1}^-(0)]_{q^{-1/2}\cdots q^{-1/2}},\right. \\
&&\qquad \left. [X_i^-(0), \cdots, X_2^-(0), X_1^-(1)]_{q^{-1/2}\cdots
q^{-1/2}}
\right]_{1\, q^{-1/2}\cdots q^{-1}\cdots q^{-1/2}\, 1}\\
&=& [X_i^-(0), \cdots, X_n^-(0), \cdots, X_{i+1}^-(0),
-zq^{n+1-i/2}\cdot \\
&&\qquad \left. \left[\Phi_{\overline {i+1}}(z),
[X_1^-(0), \cdots, X_{i-1}^-(0)]'_{q^{-1/2}\cdots q^{-1/2}}\right]_{q^{-1}}
\right]_{q^{-1/2}\cdots q^{-1}\cdots q^{-1/2}\, 1}\\
&=&zq^{(n+1)/2}\left[X_i^-(0), \cdots, X_{n-1}^-(0), \Phi_n(z),
[X_1^-(0), \cdots, X_{i-1}^-(0)]'_{q^{-1/2}\cdots q^{-1/2}}\right]_{q^{-1}\,
q^{-1/2}\cdots q^{-1/2}\, 1}\\
&=& (-1)^{n-i-1}zq^{1+i/2}\left[X_i^-(0), \Phi_{i+1}(z), [X_1^-(0), \cdots
X_{i-1}^-(0)]'_{q^{-1/2}\cdots q^{-1/2}}\right]_{q^{-1}\, 1}
\end{eqnarray*}
where we have noted that $X_{i+1}^-(0), \cdots, X_n^-(0)$ commute with
$[X_1^-(0), \cdots, X_{i-1}^-(0)]'_{q^{-1/2}\cdots q^{-1/2}}$ and
the definition of $\Phi_i(z)$'s.
\hspace{\fill} $\Box$

Now we can prove our main relation.

\noindent{\it Proof of relation \ref{E:R3}}. The idea is to move
$\Phi_{\overline 1}(z)$ to the right properly and use induction.
\begin{eqnarray*}
&& [\Phi_{\overline 1}(z), \hat{e}_0]_{q^{-1}}\\
&=& \left[[\Phi_{\overline 1}(z), X_1^-(0)]_{q^{-1/2}},
     [X_2^-(0), \cdots, X_n^-(0), \cdots, X_2^-(0),
X_1^-(1)]_{q^{-1/2}\cdots q^{-1}\cdots q^{-1/2}}\right]_{q^{-1/2}}-\\
&&\qquad q^{-1/2}\left[X_1^-(0), \cdots, X_n^-(0), \cdots, X_2^-(0),
\Phi_{\overline 2}(z)\right]_{q^{-1/2}\cdots q^{-1}\cdots
q^{-1/2}q^{1/2}}zq^n\\
&=&\left(\left[[\Phi_{\overline 1}(z), X_1^-(0), X_2^-(0)]'_{q^{-1/2}\,
q^{-1/2}}, \right. \right. \\
&&\qquad \left. [X_3^-(0), \cdots, X_n^-(0), \cdots, X_2^-(0),
X_1^-(1)]_{q^{-1/2}\cdots q^{-1}\cdots q^{-1/2}}\right]_{q^{-1/2}}+ \\
&&\qquad q^{-1/2}\left[X_2^-(0), \cdots, X_n^-(0), \cdots, X_3^-(0), \right. \\
&&\qquad \left. \left. \left[[\Phi_{\overline 1}(z), X_1^-(0)]_{q^{-1/2}},
[X_2^-(0), X_1^-(1)]_{q^{-1/2}}\right]
\right]_{q^{-1/2} \cdots q^{-1}\cdots q^{-1/2}\, 1}\right)+\\
&&\qquad (-1)^nzq^{1/2}[x_1^-(0), \Phi_2(z)]_{q^{1/2}}\\
&=&\left[[\Phi_{\overline 1}(z), X_1^-(0), X_2^-(0)]'_{q^{-1/2}\, q^{-1/2}},
\right.\\
&&\left. [X_3^-(0), \cdots, X_n^-(0), \cdots, X_2^-(0),
X_1^-(1)]_{q^{-1/2}\cdots q^{-1}\cdots q^{-1/2}}\right]_{q^{-1/2}}+\\
&& \qquad (-1)^{n-2}zq^{3/2}\left[\Phi_3(z), [X_1^-(0),
X_2^-(0)]_{q^{-1/2}}\right]_{q^{-1/2}}
\end{eqnarray*}
where we used Lemma \ref{L13} and have cancelled the term
$zq^{1/2}[x_1^-(0), \Phi_2(z)]_{q^{1/2}}$. Continuing this way for $i$ steps
we then have
\begin{eqnarray*}
&& [\Phi_{\overline 1}(z), \hat{e}_0]_{q^{-1}}\\
&=&\left[[\Phi_{\overline 1}(z), X_1^-(0), \cdots
X_{i}^-(0)]'_{q^{-1/2} \cdots q^{-1/2}},\right. \\
&& \qquad \left. [X_{i+1}^-(0), \cdots, X_n^-(0), \cdots, X_2^-(0),
X_1^-(1)]_{q^{-1/2}\cdots q^{-1}\cdots q^{-1/2}}\right]_{q^{-1/2}}\\
&&\qquad +(-1)^{n-i}z
q^{(i+1)/2}\left[\Phi_{i+1}(z), [X_1^-(0), \cdots,
X_i^-(0)]'_{q^{-1/2}\cdots q^{-1/2}}\right]_{q^{-1/2}}\\
&=&\left[[\Phi_{\overline 1}(z), X_1^-(0), \cdots
X_{n-1}^-(0)]'_{q^{-1/2} \cdots q^{-1/2}},
[X_{n}^-(0), \cdots, X_2^-(0)
X_1^-(1)]_{q^{-1/2}\cdots q^{-1/2}\, q^{-1}}\right]_{q^{-1/2}}\\
&&\qquad +(-1)^{i+1}z
q^{n/2}\left[\Phi_{n}(z), [X_1^-(0), \cdots,
X_{n-1}^-(0)]'_{q^{-1/2}\cdots q^{-1/2}}\right]_{q^{-1/2}}\\
&=&\left[[\Phi_{\overline 1}(z), X_1^-(0), \cdots
X_{n-1}^-(0)]'_{q^{-1/2} \cdots q^{-1/2}},
 q^{n/2}[X_{n}^-(1), \cdots,
X_1^-(0)]_{q^{-1/2}\cdots q^{-1/2}\, q^{-1}}\right]_{q^{-1/2}}+\\
&&\qquad -z
q^{n/2}\left[\Phi_{n}(z), [X_1^-(0), \cdots,
X_{n-1}^-(0)]'_{q^{-1/2}\cdots q^{-1/2}}\right]_{q^{-1/2}}\\
&=&q^{-n/2}\left[[\Phi_{\overline 1}(z), X_1^-(0), \cdots
X_{n-1}^-(0), X_n^-(1)]'_{q^{-1/2} \cdots q^{-1/2}\, q^{-1}},\right.\\
&&\qquad\qquad\left. [X_{n-1}^-(1), \cdots,
X_1^-(0)]_{q^{1/2}\cdots q^{1/2}}\right]_{q^{3/2}}+\\
&&\qquad q^{-n/2-1}
\left[X_n^-(1), [\Phi_{\overline 1}(z), X_1^-(0), \cdots
X_{n-1}^-(0)]'_{q^{-1/2} \cdots q^{-1/2}\, q^{-1}},\right. \\
&&\qquad \left. [X_{n-1}^-(1), \cdots, X_2^-(0),
X_1^-(0)]_{q^{1/2}\cdots q^{1/2}}\right]_{q^{1/2}\, q^2}+\\
&&-z
q^{n/2}\left[\Phi_{n}(z), [X_1^-(0), \cdots,
X_{n-1}^-(0)]'_{q^{-1/2}\cdots q^{-1/2}}\right]_{q^{-1/2}}
\quad\mbox{2nd term is $0$ by lemma \ref{L11a}}\\
&=&q^{-n/2}\left[\left[\Phi_{\overline 1}(z), [X_1^-(0), \cdots
X_{n-1}^-(0), X_n^-(1)]'_{q^{-1/2} \cdots q^{-1/2}\, q^{-1}}
\right]_{q^{-1/2}},\right. \\
&&\qquad \left. [X_{n-1}^-(0), \cdots, X_2^-(0),
X_1^-(0)]_{q^{1/2}\cdots q^{1/2}}\right]_{q^{3/2}}+\\
&&-z
q^{n/2}\left[\Phi_{n}(z), [X_1^-(0), \cdots,
X_{n-1}^-(0)]'_{q^{-1/2}\cdots q^{-1/2}}\right]_{q^{-1/2}}\\
&=&q^{-n/2}(-1)^{n-1}\left[\Phi_{\overline 1}(z), X_n^-(0), \cdots
X_{2}^-(0), X_1^-(1)]_{q^{-1/2} \cdots q^{-1/2}\, q^{-1}\, q^{-1/2}},\right.\\
&&\qquad \left. [X_{n-1}^-(0), \cdots, X_2^-(0),
X_1^-(0)]_{q^{1/2}\cdots q^{1/2}}\right]_{q^{3/2}}+\\
&&-z
q^{n/2}\left[\Phi_{n}(z), [X_1^-(0), \cdots,
X_{n-1}^-(0)]'_{q^{-1/2}\cdots q^{-1/2}}\right]_{q^{-1/2}}\\
&=&(-1)^nz\left[\Phi_{n}(z),
[X_{n-1}^-(0), \cdots, X_2^-(0),
X_1^-(0)]_{q^{1/2}\cdots q^{1/2}}\right]_{q^{3/2}}+\\
&&(-1)^{n}qz\left[\Phi_{n}(z), [X_{n-1}^-(0), \cdots,
X_{1}^-(0)]_{q^{1/2}\cdots q^{1/2}}\right]_{q^{-1/2}}\\
&=&(-1)^nz(1+q)\left[\Phi_n(z), [X_{n-1}^-(0), \cdots,
X_1^-(0)]_{q^{1/2}\cdots q^{1/2}}\right]_{q^{1/2}}\\
&=&(-1)^nz(1+q)\Phi_1(z)
\end{eqnarray*}
where we have used
$[a, b]_{q^{3/2}}+q[a, b]_{q^{-1/2}}=(1+q)[a, b]_{q^{1/2}}$.

\hspace{\fill} $\Box$

\medskip
\medskip

\centerline{Acknowledgements}
\medskip
We would like to thank Kailash Misra for conversations
on the subject. The first author is supported in part by NSA grant
MDA 904-96-1-0087. The second author is partly supported by
the Japan Society for Promotion of Science.

\end{document}